\documentclass[twocolumn,twoside]{IEEEtran}
\usepackage{amsmath,amssymb,amsthm,epsfig,dsfont,color,subfigure,empheq,graphicx,balance}
\usepackage{enumerate,url,algpseudocode,algorithm,wasysym,soul,bbm}
\usepackage[table]{xcolor}

\DeclareMathOperator{\diag}{dg}

\DeclareMathOperator{\real}{Re}
\DeclareMathOperator{\imag}{Im}

\newtheorem{remark}{Remark}

\newcommand \bzero{\mathbf{0}}
\newcommand \bone{\mathbf{1}}

\newcommand \bg{\mathbf{g}}

\newcommand \bp{\mathbf{p}}
\newcommand \bq{\mathbf{q}}

\newcommand \bs{\mathbf{s}}
\newcommand \bt{\mathbf{t}}
\newcommand \bu{\mathbf{u}}
\newcommand \bv{\mathbf{v}}
\newcommand \bw{\mathbf{w}}
\newcommand \bx{\mathbf{x}}

\newcommand \bI{\mathbf{I}}

\newcommand \bM{\mathbf{M}}

\newcommand \btheta{\boldsymbol{\theta}}

\newcommand \blambda{\boldsymbol{\lambda}}

\newcommand \bpi{\boldsymbol{\pi}}

\newcommand \bphi{\boldsymbol{\phi}}

\newcommand \mcQ{\mathcal{Q}}

\newcommand \tbu{\tilde{\mathbf{u}}}

\newcommand \cbq{\check{\mathbf{q}}}

\newcommand \bbu{\bar{\mathbf{u}}}

\begin{document}
\title{Controlling Smart Inverters using Proxies:\\
A Chance-Constrained DNN-based Approach}

\author{
Sarthak Gupta,~\IEEEmembership{Student Member,~IEEE,}
	Vassilis Kekatos,~\IEEEmembership{Senior Member,~IEEE} and
	Ming Jin,~\IEEEmembership{Member,~IEEE}

}

\mark{Gupta, Kekatos, and Jin: Controlling Smart Inverters using Proxies: A Chance-Constrained DNN-based Approach}

\maketitle

\begin{abstract}
Coordinating inverters at scale under uncertainty is the desideratum for integrating renewables in distribution grids. Unless load demands and solar generation are telemetered frequently, controlling inverters given approximate grid conditions or proxies thereof becomes a key specification. Although deep neural networks (DNNs) can learn optimal inverter schedules, guaranteeing feasibility is largely elusive. Rather than training DNNs to imitate already computed optimal power flow (OPF) solutions, this work integrates DNN-based inverter policies into the OPF. The proposed DNNs are trained through two OPF alternatives that confine voltage deviations on the average and as a convex restriction of chance constraints. The trained DNNs can be driven by partial, noisy, or proxy descriptors of the current grid conditions. This is important when OPF has to be solved for an unobservable feeder. DNN weights are trained via back-propagation and upon differentiating the AC power flow equations. An alternative gradient-free variant is also put forth, which requires only a power flow solver and avoids computing gradients. Such variant is practically relevant when calculating gradients becomes cumbersome or prone to errors. Numerical tests compare the DNN-based inverter control schemes with the optimal inverter setpoints in terms of optimality and feasibility.
\end{abstract}


\section{Introduction}\label{sec:intro}
\allowdisplaybreaks
The high penetration of DERs (such as rooftop photovoltaics, batteries, and demand response devices) introduces additional variability in distribution grid operation. Uncontrolled variations in power injections can in turn induce abrupt fluctuations in nodal voltages. Fortunately, the smart inverters interfacing DERs with the grid can propel their integration by additionally providing reactive power support. The coordinated control of DERs across a feeder can be formulated as an OPF~\cite{FCL}, \cite{ergodic}, \cite{G-SIP16}, \cite{GKS19}. If the grid is modelled using the AC power flow equations, tackling this OPF using conventional optimization-based solvers becomes formidable as DERs increase in numbers and need to be redispatched frequently under dynamic conditions. Therefore, operators may not have the time and computational resources to solve an AC-OPF every few seconds. At the same time, solving an OPF presumes all problem inputs (load demands and solar generation) to be precisely known. Nonetheless, such parameters are oftentimes described stochastically, observed under noise and delays, or the operator can monitor only few of them in real time. Therefore, even if an operator can afford solving an AC-OPF every few minutes, it may not know all loads and solar generation at the level required by the AC-OPF. The need to compute reasonable DER control decisions in real time and without knowing the current grid conditions in full detail is the driving motivation of this work.

Alternatively, recent literature advocates the use of machine learning (ML) models to solve minimization problems under the \emph{learning-to-optimize} paradigm. Due to the rich modeling and fast inference capabilities of ML models, learning-to-optimize becomes relevant to scenarios where large-scale non-linear optimization tasks have to be solved frequently enough and/or under uncertain or partially observed inputs. ML-based schemes for tackling the OPF have already been explored and can be broadly classified into the \emph{OPF-then-learn} and \emph{OPF-and-learn} categories. Methods within the former category involve two steps. They first generate a labelled training dataset by solving a large number of OPFs. The features and labels of this dataset consist of the OPF input parameters and outputs (optimal DER setpoints), respectively. During the second step, an ML model is trained to predict the dataset labels in the conventional supervised manner. Within the \emph{OPF-then-learn} category, kernel--based regression has been employed to learn inverter control rules in~\cite{Kara18}, physics-informed stacked extreme learning has been utilized for OPFs~\cite{SELM21}, while DNNs have been trained to predict OPF solutions under a linearized~\cite{DeepOPFPan19,pascal21} and the exact AC grid models~\cite{ZamzamBaker19}, \cite{GuhaACOPF}, \cite{ZhangDNNOPF}, \cite{RibeiroGNNOPF}.To expedite the first step, the sensitivity-informed learning method of \cite{SGKCB2020} and \cite{L2O2021} trains a DNN to match not only the OPF minimizers, but also their partial derivatives with respect to the OPF inputs. Despite the data efficiency enhancement offered by sensitivity-informed learning, the OPF-then-learn strategy lacks feasibility guarantees and presumes OPF input parameters are deterministic and known. Furthermore, generating labels by solving several OPFs incurs a significant computational overhead. Consequently, the OPF-then-learn strategy is not well suited for scenarios where the optimal policies need to be re-learned frequently on account of changing underlying data distribution.

{Instead of this two-step approach of first generating OPF labels and then training the ML model to fit these labels in a least-squares (LS) sense, \emph{OPF-and-learn} approaches learn the ML model directly while solving the OPF in a single step. This is achieved by altering the optimization algorithm used for training the ML model. Rather than fitting OPF solutions in the LS sense, the training algorithm uses the objective or the Lagrangian function of the OPF at hand. By training the ML model directly using the cost and constraint functions of a stochastic OPF, we avoid the computationally expensive step of solving numerous OPF instances beforehand to generate a labeled dataset. Because the ML model is trained now over several OPF instances, it is particularly suited for \emph{stochastic} OPF formulations, where one is interested on the average or probabilistic performance of the learned OPF decisions or policies over uncertain and/or time-varying conditions.} Hence, methods within the \emph{OPF-and-learn} category are more appropriate for dynamic applications where ML models need to be continuously retrained. Under the \emph{OPF-and-learn} paradigm, reference~\cite{JKGD19} adopts support vector machines (SVMs) to design inverter control rules, adjusted to grid conditions in a quasi-stationary fashion. Albeit SVM-based rules can be learned via a convex program, kernel functions have to be specified beforehand. In~\cite{Yang19}, inverter control rules are optimized along with capacitor status decisions to minimize voltage deviations using a two-timescale reinforcement learning (RL) approach yet no feeder-level constraints are involved. Enforcing network constraints is challenging for ML-based OPF methods. One could heuristically project the ML prediction for the OPF solution~\cite{ZamzamBaker19}, \cite{JKGD19}, or consider penalizing constraint deviations~\cite{Kara18}, \cite{DeepOPFPan19}, \cite{Yang19}, or include deterministic constraints per training sample and solve using dual approximation~\cite{pascal21}.

An alternative for dealing with constraints in the learning-to-optimize process is through the discounted return functions used in RL approaches. In this context, reference~\cite{ZhangISU} updates DNN-based inverter policies continuously by successively linearizing feeder constraints. A similar \emph{safe RL} learning scheme is put forth in~\cite{Nanpeng19}, which focuses on regulating the number of voltage violations across nodes through a method of multipliers strategy. However, both these works updated the policy parameters by solving an optimization problem at every update step that might be computationally intensive, restricting their applicability on dynamic scenarios. Secondly, both~\cite{ZhangISU} and \cite{Nanpeng19} focus on scenarios where measurements across feeder nodes are available to policies in real time. Finally, RL-based approaches are in general more complex in implementation, which can be hindering their adoption by grid operators; e.g., the \emph{safe RL} strategy of \cite{Nanpeng19} involves nine different DNNs.

Primal-dual learning~\cite{Ribeiro19}, provides a computationally less intensive alternative to RL for handling feeder constraints. Different from~\cite{ZhangISU} and \cite{Nanpeng19}, primal-dual learning involves simpler gradient-based updates of the policy parameters and the related dual variables alike. Stochastic primal-dual updates were first applied towards learning the optimal control policies for smart inverters to enforce averaged voltage constraints in the conference precursor of our work~\cite{SG20}. 

\emph{Contributions:} Building on the \emph{OPF-and-learn} approach of \cite{SG20}, this work trains a DNN that learns a stochastic inverter control policy to provide near-optimal setpoints in real time and without fully knowing the current grid conditions. The contributions over \cite{SG20} extend in four fronts. Firstly, the underlying feeder is represented using the exact AC model rather than the approximate linearized model~\cite{TJKT20} previously employed in~\cite{SG20}. This extension is non-trivial as the gradients needed for the stochastic primal-dual updates of the DNN are now found in an indirect manner using the underlying AC power flow equations and the inverse function theorem. Secondly, we illustrate the versatility of stochastic primal-dual updates as they can cope with probabilistic voltage constraints via convex restrictions~\cite{Nemirovski07}. Thirdly, we demonstrate how primal-dual learning can also be used when the DNN has to be fed with partial information during operation. This adheres to practical setups where the utility might have real-time telemetry only over a subset of grid locations. Fourthly, we propose gradient-free counterparts of the primal-dual updates that train the DNN knowing only the values of voltages and losses, and not their gradients with respect to the control variables. Such approaches are useful when the feeder model is known but complex (due to the presence of transformers, capacitors, ZIP loads) and differentiation becomes perplex, but the utility has access to a power flow solver.

\emph{Paper outline:} The rest of the paper is organized as follows. Section~\ref{sec:opf} formulates the task of DNN-based smart inverter control, and puts forth an \emph{averaged} and a \emph{probabilistic} scheme. Section~\ref{sec:train} elaborates on primal-dual DNN learning for both schemes. Section~\ref{sec:train} calculates the gradients needed for the stochastic updates, while the gradient-free implementation is presented in Section~\ref{sec:modelfree}. The novel DNN-based inverter control strategies are evaluated using real-world data on the IEEE 37-bus feeder in Section~\ref{sec:tests}. Conclusions are drawn in Section~\ref{sec:conclusions}.

\emph{Notation:} Lower- (upper-) case boldface letters denote column vectors (matrices), and calligraphic symbols are reserved for sets. Symbol $^{\top}$ stands for transposition and $\|\bx\|_2$ denotes the $\ell_2$-norm of $\bx$. Vectors $\bzero$ and $\bone$ are respectively the vectors of all zeros and ones of appropriate dimensions.

\section{Problem Formulation}\label{sec:opf}
Consider a feeder with $N+1$ buses. The substation is indexed by $0$ and the remaining buses comprise the set $\mathcal{N}:=\{1,\dots,N\}$. 
Let $p_n+jq_n$ be the complex power injection at bus $n$. Its active power component can be decomposed as $p_n=p_n^g-p_n^c$, where $p_n^g$ is the solar generation and $p_n^c$ the inelastic load at bus $n$. Its reactive power component can be similarly expressed as $q_n=q_n^g-q_n^c$. The vectors $(\bp,\bq)$ collecting the power injections at all non-substation buses can be decomposed as $\bp=\bp^g-\bp^c$ and $\bq=\bq^g-\bq^c$. For simplicity, it is assumed that each bus hosts at most one inverter, and so index $n$ will be henceforth used for buses and inverters interchangeably.

Let vector parameter $\btheta:=\{\bp^c,\bq^c,\bp^g\}\subseteq \mathbb{R}^{3N}$ collect the loads (active and reactive) and active solar generation at all non-substation buses. We will henceforth term $\btheta$ as the vector of \emph{grid conditions}.  Given $\btheta$, the task of reactive power control by DERs aims at optimally setting $\bq^g$ to minimize a feeder-wide objective while complying with network and inverter limitations. Starting with the latter, the reactive power injected by inverter $n$ is limited by its given apparent power limit $\bar{s}_{n}$. Apparent power constraints are local and will be collectively denoted by 
\begin{equation}\label{eq:Q}
\bq^g\in\mathcal{Q}_{\btheta}:=\left\{\bq^g:|q_{n}^g|\leq  \sqrt{\overline{s}_{n}^2-(p_{n}^g)^2}~\forall n\right\}.
\end{equation}
where the subscript in $\mathcal{Q}_{\btheta}$ denotes that the feasible space changes as the value of solar generation $\bp_g$ changes.

The task of coordinating inverters can be centrally handled by the utility operator. The operator finds the reactive power setpoints for DERs by minimizing ohmic losses on distribution lines while maintaining voltage magnitudes within per-bus bounds $[\underline{\bv},\overline{\bv}]$ as
\begin{align}\label{eq:opf}
\min_{\bq\in\mcQ_{\btheta}}~&~\ell(\bq,\btheta)\\
\mathrm{s.to}~&~\underline{\bv}\leq\bv(\bq,\btheta)\leq \overline{\bv}.\nonumber
\end{align}
 where $\bv$ is the vector of bus voltage magnitudes. We will henceforth refer to voltage magnitudes as voltages unless stated otherwise. We slightly abuse notation and use $\bq$ in lieu of $\bq^g$ to unclutter notation, as $\bq^c$ is included in $\btheta$ anyway. Note that expressions $\ell(\bq,\btheta)$ and $\bv(\bq,\btheta)$ capture the dependence of losses and nodal voltages on the reactive setpoints of DERs $\bq$ under the current grid conditions $\btheta$. {It is assumed that the feeder model and the participating inverters are known and remain fixed throughout the control period.}

Solving \eqref{eq:opf} can be computationally and communication-wise taxing if $\btheta$ changes frequently. Moreover, by the time \eqref{eq:opf} is solved and optimal setpoints are downloaded to DERs, grid conditions $\btheta$ may have changed rendering the computed setpoints obsolete~\cite{KeWaCoGia15}, \cite{ergodic}. To account for the uncertainty in $\btheta$, we propose two stochastic formulations. The first formulation replaces $\ell(\bq,\btheta)$ and $\bv(\bq,\btheta)$ with their averages:
\begin{align}\label{eq:opf2}
\min_{\bq\in\mcQ_{\btheta}}~&~\mathbb{E}[\ell(\bq,\btheta)]\\ 
\mathrm{s.to}~&~\underline{\bv}\leq\mathbb{E}[\bv(\bq,\btheta)]\leq \overline{\bv}\nonumber
\end{align}
where the expectation $\mathbb{E}$ is with respect to $\btheta$. We refer to \eqref{eq:opf2} as the \emph{averaged formulation}. While the averaged formulation takes care of uncertainties in $\btheta$, the obtained setpoints may violate the voltage limits in \eqref{eq:opf} quite frequently. This undesirable behavior results from the fact that constraining the average value of voltages does not provide strong guarantees on their per-instance values. Nevertheless, the averaged formulation has an attractive structure that permits straightforward stochastic gradient descent (SGD)-based steps to arrive at the optimal setpoints as we will see later.

A more conservative approach is possible through the \emph{probabilistic formulation}
\begin{subequations}\label{eq:opf3}
\begin{align}
\min_{\bq\in\mcQ_{\btheta}}~&~\mathbb{E}[\ell(\bq,\btheta)]\\ 
\mathrm{s.to}~&~\Pr\left[\underline{v}_n\geq v_n(\bq,\btheta)\right]\leq \alpha,~\forall n\in\mathcal{N}\label{eq:chancelow}\\
~&~\Pr\left[\overline{v}_n\leq v_n(\bq,\btheta) \right]\leq \alpha,~\forall n\in\mathcal{N}\label{eq:chancehigh}
\end{align}
\end{subequations}
where \eqref{eq:chancelow}--\eqref{eq:chancehigh} ensure each bus voltage remains within the desired limits with a probability of at least  $1-\alpha$ on each side. Here $\alpha\in(0,1)$ is a small \emph{violation probability}. In contrast to \eqref{eq:opf2}, the formulation in \eqref{eq:opf3} focuses on restricting the frequency of occurrence of voltage violations. Problem \eqref{eq:opf3} can be rewritten as
\begin{subequations}\label{eq:opf4}
\begin{align}
\min_{\bq\in\mcQ_{\btheta}}~&~\mathbb{E}[\ell(\bq,\btheta)]\\ 
\mathrm{s.to}~&~\mathbb{E}\left[\mathbbm{1}(\underline{v}_n- v_n(\bq,\btheta))\right]\leq \alpha,~\forall n\in\mathcal{N}\label{eq:chancelow4}\\
~&~\mathbb{E}\left[\mathbbm{1}(v_n(\bq,\btheta)- \overline{v}_n)\right]\leq \alpha,~\forall n\in\mathcal{N}\label{eq:chancehigh4}
\end{align}
\end{subequations}
where the indicator function $\mathbbm{1}(x)$ is defined as 
\begin{equation}\label{eq:indicator}
\mathbbm{1}(x)=\left\{\begin{array}{ll}
1&,~x\geq 0\\
0&,~x<0
\end{array}
\right..
\end{equation}
The probabilistic formulation is difficult to handle since the indicator function is neither convex nor differentiable. In quest of workable alternatives, Section~\ref{sec:train} pursues convex approximations of constraints \eqref{eq:chancelow4}--\eqref{eq:chancehigh4}. 

For now, let both formulations be represented by the general stochastic program
\begin{align}\label{eq:opf5}
\min_{\bq\in\mcQ_{\btheta}}~&~\mathbb{E}[\ell(\bq,\btheta)]\\ 
\mathrm{s.to}~&~\mathbb{E}[\bg(\bq,\btheta)]\leq \bzero.\nonumber
\end{align}
Note that solving \eqref{eq:opf5} results in a single `one-size-fits-all' $\bq$ that does not adapt to different $\btheta$'s. To render DER setpoints responsive to grid conditions, we resort to a \emph{control policy}, where the reactive setpoints $\bq$ are captured by a function $\bq=\pi(\btheta;\bw)$, which is parameterized by $\bw$.


Ideally, the control policy is driven by the vector of grid conditions $\btheta$. Nevertheless, during real-time operation, the operator controlling the DERs may not be able to observe the complete $\btheta$. Instead, it may have to act upon a proxy $\bphi$ of the actual $\btheta$. The DER control policy driven by $\bphi$ can then be found by solving the constrained stochastic minimization
\begin{align}\label{eq:opf6}
\min_{\bw:\bpi(\bphi;\bw)\in\mcQ_{\btheta}}~&~\mathbb{E}[\ell(\bpi(\bphi;\bw),\btheta)]\\
\mathrm{s.to}~&~\mathbb{E}[\bg(\bpi(\bphi;\bw),\btheta)]\leq \bzero.\nonumber
\end{align}
The DER control policies found through \eqref{eq:opf6} are adaptive to the proxy vector $\bphi$ and the optimization is over the parameters $\bw$. Policies account for the uncertainty over $\btheta$, and correspondingly $\bphi$. Note that the expectations in \eqref{eq:opf6} couple the system's performance across OPF instances of $\btheta$. The notation $\ell(\bpi(\bphi;\bw),\btheta)$ captures the fact that the control policy is fed by proxy $\bphi$ to determine $\bq$, but of course ohmic losses depend on the actual grid conditions $\btheta$.

The proxy vector $\bphi$ can be chosen to represent the operational setup for which the control policies are being designed. In the absence of real-time measurements from all nodes, and/or to save on communication overhead, vector $\bphi$ can consist of active line flows from distribution lines~\cite{SG20}. Meteorological data such as solar irradiance and ambient temperature, which serve as surrogates for $\bp$, can also be included in $\bphi$. One can also explore convolutional neural networks (CNNs)-based policies that accept sky images in place of solar irradiance measurements as inputs to be included in $\bphi$. {Similarly, the proxy vector $\bphi$ can also represent partial, delayed, or noisy data on the grid conditions, or even aggregate versions of them.}  In Section~\ref{sec:tests}, a more straightforward scenario is explored whereby measurements from a subset of buses in $\mathcal{N}$ are assumed to be available in real-time, resulting in $\bphi\subset \btheta$.

Previous works have studied linear inverter control policies of the form $\bpi(\bphi;\bw)=\bw^\top\bphi$; see e.g.,~\cite{Jabr18}, \cite{LinBitar18}, \cite{Baker18}. Nonetheless, the optimal policy $\bpi(\bphi;\bw)$ is not necessarily affine in $\bphi$, especially when $\bphi$ is a proxy for $\btheta$. The grand challenge towards scalable inverter control is to design \emph{nonlinear control policies}. To this end, in~\cite{JKGD19}, we modeled $\bpi(\bphi;\bw)$ as a support vector machine (SVM) and designed the policy through an OPF formulation. The advantage of SVM-based policies is that they can be trained to optimality using convex optimization. Nonetheless, selecting the appropriate kernel and control input $\bphi$ can be challenging. Inspired by their field-changing performance in various engineering tasks, here we propose modeling the DER control policy $\bq=\bpi(\bphi;\bw)$ by a DNN. The proxy vector $\bphi$ of grid conditions $\btheta$ is fed as an input to the DNN. Vector $\bw$ carries the weights of the DNN across all layers. The output of the DNN $\bpi(\bphi;\bw)$ predicts the sought inverter setpoints $\bq$. Figure~\ref{fig:activation} provides a schematic of the architecture. We propose learning weights $\bw$ in a data-driven physics-aware fashion.

\begin{figure}[t]
    \centering
    \includegraphics[scale=0.35]{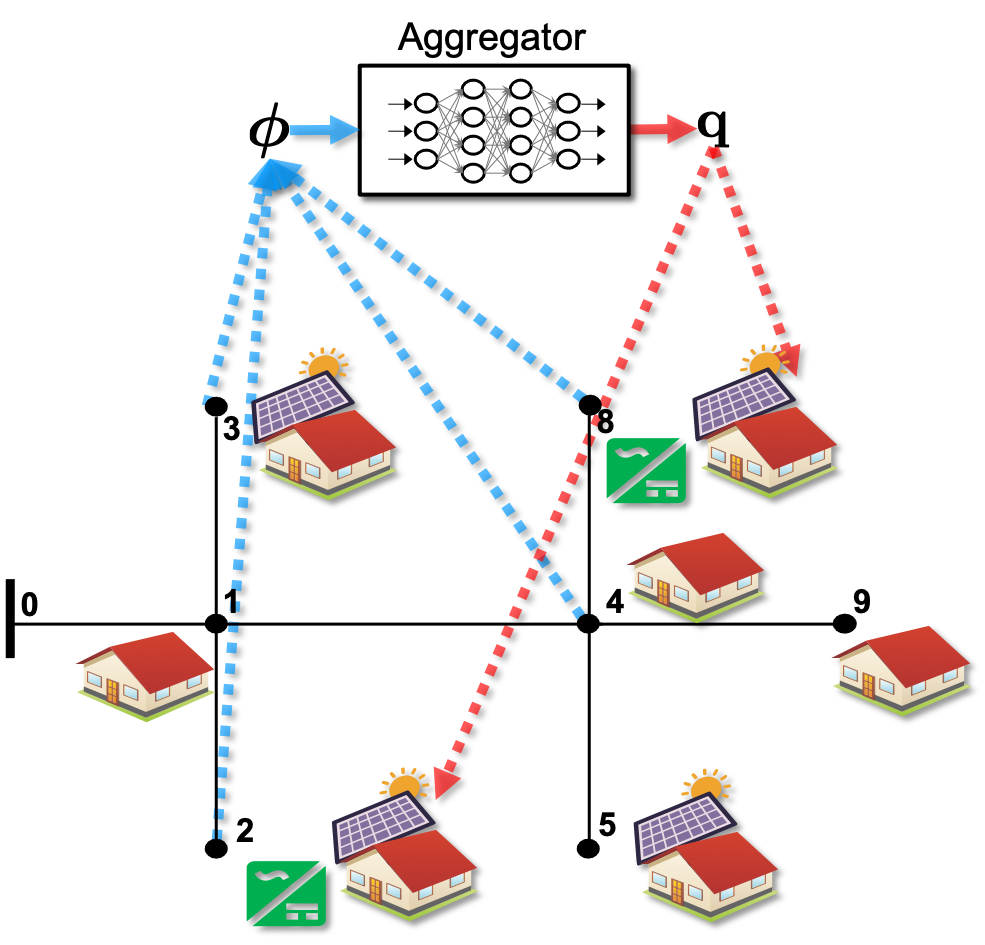}
    \caption{Real-time operation of the proposed DER inverter control scheme. Proxy vector $\bphi$ consisting of measurements from nodes $\{2,3,4,8\}$ is transmitted to the operator. The DNN-based policy acts upon $\bphi$ to predict setpoints $\bq$, which are then broadcast to the inverters at buses 2 and 8.}
    \label{fig:operation}
\end{figure}

Figure~\ref{fig:operation} depicts the real-time operation of the proposed control strategy. The smart inverters to be controlled are located on buses $2$ and $8$. The proxy vector $\bphi\subset\btheta$, consisting of $\{p^c_n,q^c_n,p^g_n\}$ measurements from buses 2, 3, 4 and 8, is transmitted to the operator. The operator feeds $\bphi$ as an input to the DNN-based policy and predicts setpoints $\bq$. These setpoints are then broadcast to DER inverters for implementation. It is worth emphasizing here that although the operator needs to know pairs of $(\btheta,\bphi)$ during training, the DNN-based policy operates solely on $\bphi$.


\section{Primal-Dual DNN Training }\label{sec:train}
The stochastic formulation in \eqref{eq:opf6} is challenging to solve on account of the expectation operator in both the objective and constraints. Computing the needed expectations requires knowing the probability density functions (pdf) of $\bphi$ and $\btheta$. Even if these pdfs are known, computing the expectations is still non-trivial granted the policies $\bpi(\bphi;\bw)$ are non-linear in $\bphi$. These complications promulgate a stochastic approximation approach towards solving \eqref{eq:opf6}. In the conventional machine learning setup, the weights of a DNN are found by minimizing a data-fitting loss function under no constraints via stochastic gradient descent. Here, to accommodate constraints, we adopt the stochastic primal-dual updates of~\cite{Ribeiro19} as presented next.

Consider the Lagrangian function of the problem in \eqref{eq:opf6}
\begin{align}\label{eq:lagrangian}
L(\bw;\blambda):=  \mathbb{E}[\ell(\bpi(\bphi;\bw),\btheta)] + \blambda^\top\mathbb{E}[\bg(\bpi(\bphi;\bw),\btheta)]
\end{align} where $\blambda$ is the vector of Lagrange multipliers corresponding to the constraints in \eqref{eq:opf6}. Vector $\blambda$ concatenates the multipliers $\overline{\blambda}$ and $\underline{\blambda}$ associated with the lower and upper voltage limits in \eqref{eq:opf2} and \eqref{eq:opf4}. A stationary point for the related dual problem
\begin{align}\label{eq:dual}
D^*:=\max_{\blambda\geq \bzero}\min_{\bw:\bpi(\bphi;\bw)\in\mcQ_{\btheta}} L(\bw;\blambda)
\end{align} 
can be obtained iteratively using the primal-dual updates indexed by $k$~(cf.~\cite{Ribeiro19}): 
\begin{subequations}\label{eq:pdupdate}
\begin{align}
    \bw^{k+1}&:=\big[\bw^{k}-\mu_w\nabla_{\bw}L(\bw^{k};\blambda^{k})\big]_{\mcQ_{\btheta}}\label{eq:pdupdate:p}\\
    \blambda^{k+1}&:=\big[\blambda^{k}+\mu_{\lambda}\nabla_{\blambda}L(\bw^{k+1};\blambda^{k})\big]_+\label{eq:pdupdate:d}
\end{align}
\end{subequations}
where $(\mu_w,\mu_{\lambda})$ are positive step sizes. Here primal variables are updated through projected gradient descent steps on the Lagrangian function. Dual variables are updated through projected gradient ascent steps again on the Lagrangian function. The operator $[x]_+=\max\{x,0\}$ is applied entry-wise and ensures $\blambda\geq\mathbf{0}$ at all times. 
\begin{figure}[t]
    \centering
    \includegraphics[scale=0.4]{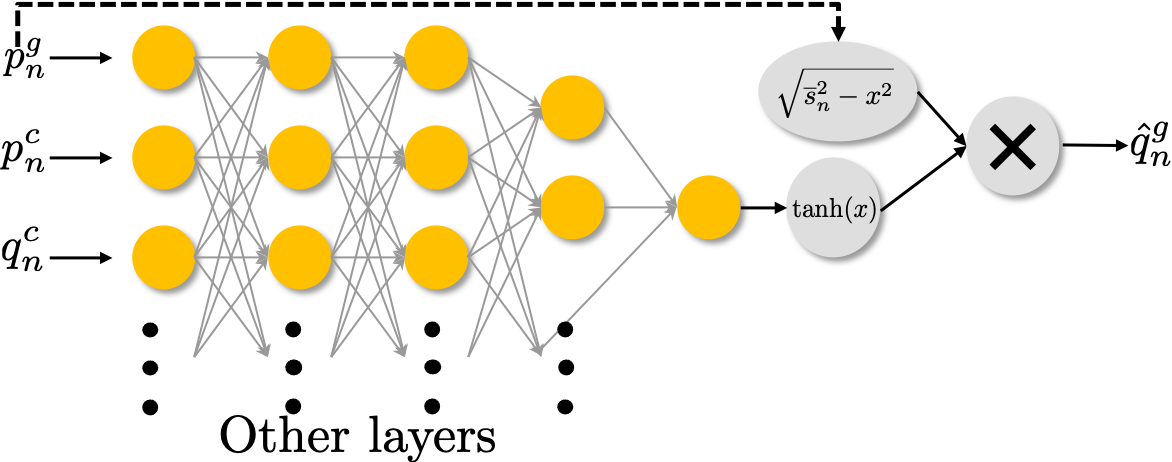}
    \caption{The inverter control policy $\bq=\bpi(\bphi;\bw)\in\mcQ_{\btheta}$ has been implemented using a DNN. Vertically stacked nodes represent the nonlinear activation functions of a single layer. Edges across nodes correspond to weights applied linearly on node output to compute the inputs to the next layer. The vector of grid conditions $\btheta$ (or its proxy $\bphi$) is fed as the input to the DNN. Vector $\bw$ collects all weights of the DNN. The DNN output provides the inverter setpoints $\bq$. The activation function of the output layer of the DNN has been modified to ensure $\bpi(\bphi;\bw)\in\mcQ_{\btheta}$ for all $\bphi$. The output neuron feeds into the $\tanh()$ activation function and then scaled by $\sqrt{\overline{s}_{n}^2-(p_{n}^g)^2}$. The scaling operation involves a skip connection from $p_n^g$, which is one of the inputs to the DNN.}
    \label{fig:activation}
\end{figure}

The operator $[\cdot]_{\mcQ_{\btheta}}$ projects $\bw^{k+1}$ such that $\bpi(\bphi;\bw^{k+1})\in\mcQ_{\btheta}$ for all $\bphi$. A direct way to confine the DNN output $q_n^g$ within $\pm \sqrt{\overline{s}_{n}^2-(p_{n}^g)^2}$ is to use the hyperbolic tangent ($\tanh$) as the output activation function and then scale the output by $ \sqrt{\overline{s}_{n}^2-(p_{n}^g)^2}$. While the apparent power limit $\overline{s}_{n}$ is known \emph{a priori}, the solar generation $p_{n}$ is available to the DNN as an input. The required scaling is easily accommodated by minor architectural modifications to the activation layer of the DNN. Figure~\ref{fig:activation} illustrates this process for a single-output neuron. Ensuring that $\bpi(\bphi;\bw)\in\mcQ_{\btheta}$ at all times obviates the need for projecting the weight updates henceforth. Note that the gradient $\nabla_{\blambda}L(\bw;\blambda)$ in the dual variable update in \eqref{eq:pdupdate:d} can be substituted as $\nabla_{\blambda}L(\bw;\blambda)=\bg(\bpi(\bphi;\bw),\btheta)$.


Following a stochastic approximation approach, the ensemble averages in \eqref{eq:lagrangian} are first surrogated by sample averages computed over a set of $S$ scenarios $\{\bphi^s,\btheta^s\}_{s=1}^S$. The average ohmic losses for example can be approximated as 
\begin{align}\label{eq:sample_average}
\mathbb{E}[\ell(\bpi(\bphi;\bw),\btheta)]\simeq \frac{1}{S}\sum_{s=1}^S\ell(\bpi(\bphi^s;\bw),\btheta^s).
\end{align}
Even with the sample approximation in \eqref{eq:sample_average}, computing the gradients needed in \eqref{eq:pdupdate} remains computationally expensive as one needs to compute gradients for each one of the $S$ training examples. Taking ohmic losses for example, we have that
\begin{equation}\label{eq:SGD1}
\mathbb{E}[\nabla_{\bw} \ell(\bpi(\bphi;\bw^k),\btheta)]\simeq  \frac{1}{S}\sum_{s=1}^S \nabla_{\bw}\ell(\bpi(\bphi^s;\bw^k),\btheta^s).
\end{equation}
Notice the two indices in~\eqref{eq:SGD1}: index $k$ indexes primal/dual updates, and index $s$ indexes data (scenarios). To perform iteration $k$, one has to cycle across all scenarios $s=1,\ldots,S$, and compute the derivatives of losses with respect to the previous update $\bw^k$ for each one of the scenarios. Stochastic approximation alleviates this burden by approximating the gradients needed in \eqref{eq:pdupdate} using a \emph{single} scenario per iteration. In other words, gradients are approximated not by \eqref{eq:SGD1}, but using a single datum as
\begin{equation}\label{eq:SGD2}
\mathbb{E}[\nabla_{\bw} \ell(\bpi(\bphi;\bw^k),\btheta)]\simeq  \nabla_{\bw}\ell(\bpi(\bphi^s;\bw^k),\btheta^s).
\end{equation}
Therefore, at iteration $k$, stochastic approximation selects a scenario $s$ to compute all gradients needed in \eqref{eq:pdupdate}. Scenarios can be selected at random or sequentially. Either way, since each iteration $k$ ends up using only a single scenario $s$, we will henceforth use symbol $k$ to index both iterations and scenarios. This is without loss of generality as $(\bphi^k,\btheta^k)$ simply means the scenario picked at iteration $k$. 

\begin{figure*}[t]
    \centering
    \includegraphics[scale=0.7]{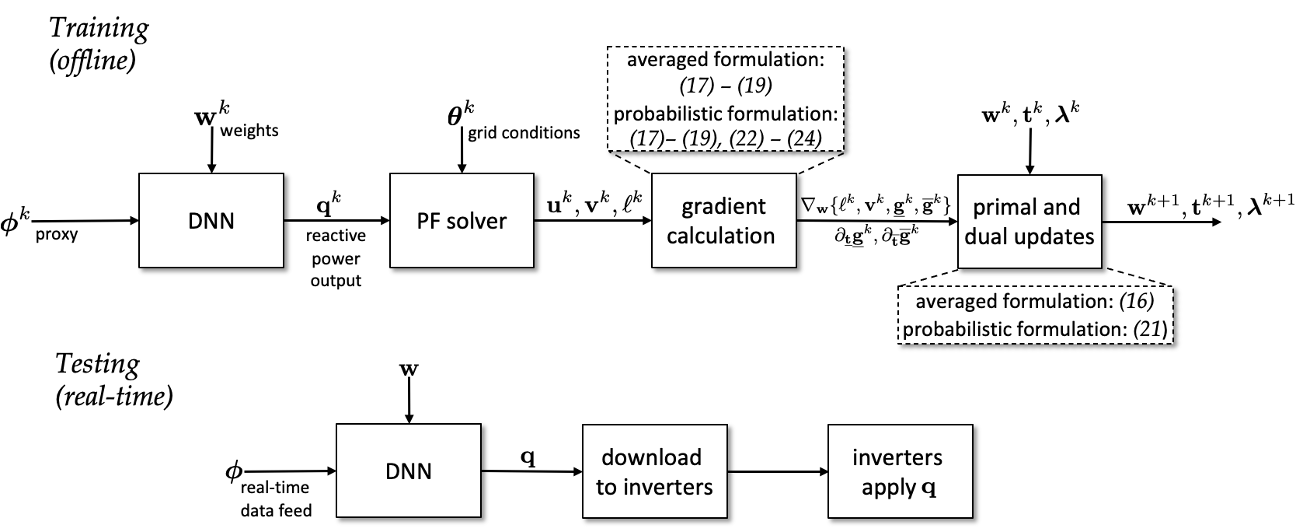}
    \caption{Workflow for the training and testing (operation) phases of the proposed DNN-based inverter control strategy.}
    \label{fig:train_vs_test}
\end{figure*}

Given the aforesaid simplification, the gradients in \eqref{eq:pdupdate} can be approximated using a single scenario per update as~\cite{Ribeiro19}
\begin{subequations}\label{eq:spdupdate}
\begin{align}
    \bw^{k+1}&:=\bw^{k}-\mu_w\left(\nabla_{\bw}\ell^k+ \big(\nabla_{\bw}{\bg^k}\big)^\top\blambda^k\right)\label{eq:spdupdate:p}\\
    \blambda^{k+1}&:=\left[\blambda^{k}+\mu_{\lambda} \bg\left(\bpi(\bphi^k;\bw^{k+1}),\btheta^k\right)\right]_+\label{eq:spdupdate:d}
\end{align}
\end{subequations}
where the shorthand notation $\nabla_{\bw}\ell^k$ denotes the gradient of $\ell$ and $\nabla_{\bw}{\bg^k}$ the Jacobian matrix of $\bg$, both with respect to $\bw$ and both evaluated at $(\bphi^k,\bw^k,\btheta^k)$. The rest of this section explains how the gradients appearing in \eqref{eq:spdupdate:p} can be computed for the averaged and probabilistic formulations, while Figure~\ref{fig:train_vs_test} summarizes the workflow for the training and testing (operational) phases for both formulations.

\subsection{Averaged Formulation}\label{subsec:ave}
For the averaged formulation, the stochastic primal-dual updates can be obtained by replacing $\bg(\bpi(\bphi;\bw),\btheta)$ with the constraint functions from \eqref{eq:opf2} to get
\begin{subequations}\label{eq:avg}
\begin{align}
    \bw^{k+1}&:=\bw^{k}-\mu_w\left(\nabla_{\bw}\ell^k+\left(\nabla_{\bw}\bv^k\right)^\top(\overline{\blambda}^k-\underline{\blambda}^k)\right)\label{eq:avg:p}\\
     \underline{\blambda}^{k+1}&:=\left[\underline{\blambda}^{k}+\mu_{{\lambda}}\left(\underline{\bv}-\bv\left(\bpi(\bphi^k;\bw^{k+1}),\btheta^k\right)\right)\right]_+\label{eq:avg:dup}\\
    \overline{\blambda}^{k+1}&:=\left[\overline{\blambda}^{k}+\mu_{{\lambda}}\left(\bv\left(\bpi(\bphi^k;\bw^{k+1}),\btheta^k\right)-\overline{\bv}\right)\right]_+.\label{eq:avg:dlow}
\end{align}
\end{subequations}
To compute the Jacobian matrix $\nabla_{\bw}\bv$ of voltages with respect to DNN weights, we resort to the power flow equations. Let vector $\tbu\in\mathbb{C}^{N+1}$ collect the complex voltage phasors at all buses and define $\bbu:=[\real(\tbu)^\top~\imag(\tbu)^\top]^\top\in\mathbb{R}^{2N+2}$. Vector $\bu\in\mathbb{R}^{2N}$ is obtained by dropping the first and $(N+1)$-th entries of $\bbu$. These two entries correspond to the substation voltage $\tilde{u}_0$, which is assumed constant. The sought Jacobian matrix can be computed via the chain rule as
\begin{equation}\label{eq:jacobvolt}
\nabla_{\bw}{\bv}=\nabla_{\bq}{\bv}\cdot\nabla_{\bw}{\bq}=\nabla_{\bu}{\bv}\cdot\nabla_{\bq}{\bu}\cdot\nabla_{\bw}{\bq}.
\end{equation}

We next elaborate on the three Jacobian matrices needed in \eqref{eq:jacobvolt}. The $N\times 2N$ matrix $\nabla_{\bu}{\bv}$ can be readily computed by definition of voltage magnitudes. Its $(n,m)$-th entry is
\[[\nabla_{\bu}{\bv}]_{n,m}=\left\{\begin{array}{ll}
\frac{u_n}{v_n} &,~m=n\\
\frac{u_{n+N}}{v_n} &,~m=n+N\\
0&,~\text{otherwise}.
\end{array}
\right.\]

We proceed with finding $\nabla_{\bq}{\bu}$. If $\bp+j\bq$ is the vector of complex power injections at all buses excluding the substation, define $\bs:=[\bp^\top~\bq^\top]^\top$. Per the power flow equations, every entry $i$ of $\bs$ can be expressed as a quadratic function of $\bbu$, that is $s_i=\bbu^\top\bM_i\bbu$ for a symmetric real-valued matrix $\bM_i$ derived from the bus admittance matrix of the feeder; see e.g.,~\cite{PSSE-Redux} for the detailed expressions for $\bM_i$'s. Therefore, the $i$-th row of the Jacobian matrix $\nabla_{\bbu}\bs$ is given by $2\bbu^\top\bM_i$. By dropping the first and $(N+1)$-th column of $\nabla_{\bbu}\bs$, we obtain $\nabla_{\bu}\bs$. Under mild technical conditions, the inverse function theorem predicates that
\begin{equation}\label{eq:ift}
\nabla_{\bs}{\bu}=(\nabla_{\bu}\bs)^{-1}.
\end{equation}
Matrix $\nabla_{\bq}{\bu}$ can be clearly obtained by keeping only the last $N$ rows of $\nabla_{\bs}{\bu}$. The third matrix $\nabla_{\bw}\bq$ can be readily computed using gradient back-propagation across the DNN.

The gradient vector of ohmic losses with respect to the DNN weights can be computed similarly as
\begin{equation}\label{eq:lossgrad}
(\nabla_{\bw}\ell)^\top=(\nabla_{\bu}\ell)^\top \cdot \nabla_{\bq}{\bu}\cdot\nabla_{\bw}{\bq}.
\end{equation}
The gradient $\nabla_{\bu}\ell$ can be easily computed by recognizing
\begin{align*}
   \ell=\sum_{n=0}^Np_n=\bbu^\top\left( \sum_{n=0}^N\bM_i\right)\bbu
\end{align*}

where matrix $\bM_0$ describes the power injection at the substation as $p_0=\bbu^\top\bM_0\bbu$ similar to the remaining injections. It is then obvious that $\nabla_{\bbu}\ell=2\sum_{n=0}^N\bM_i \bbu$. The gradient $\nabla_{\bu}\ell$ of \eqref{eq:lossgrad} is found by dropping the first and $(N+1)$-th entries of vector $\nabla_{\bbu}\ell$.


\subsection{Probabilistic Formulation}
For the probabilistic formulation of \eqref{eq:opf4}, the update steps in \eqref{eq:spdupdate} cannot be applied directly. This is because constraints \eqref{eq:chancelow4}--\eqref{eq:chancehigh4} involve the indicator function that is not differentiable, and so the Jacobian matrix $\nabla_\bw\bg$ does not exist. Moreover, these constraints are non-convex, thus prohibiting stochastic (sub)gradient updates. To circumvent these complications, we instead turn to the  convex CVaR approximations to \eqref{eq:opf4} using stochastic subgradient primal-dual updates.



We first briefly review the CVaR approximation of chance constraints from~\cite{Nemirovski07}, and then compute the related subgradients. For some $\alpha\in(0,1)$, consider the chance constraint $\Pr[f_\theta(x)\geq 0]=\mathbb{E}_{\theta}[\mathbbm{1}(f_\theta(x))]\leq \alpha$, where function $f$ depends on the optimization variable $x$ and a random variable $\theta$. Note that $\mathbbm{1}(f_\theta(x))\leq [1+f_\theta(x)/t]_+$ for all $f_\theta(x)$ and $t>0$. This is easy to verify by checking the two cases in the definition of the indicator function in~\eqref{eq:indicator}. Then, if there exists a $t>0$ satisfying $\mathbb{E}_\theta\left[[1+f_\theta(x)/t]_+\right]\leq \alpha$, the original chance constraint $\mathbb{E}_{\theta}[\mathbbm{1}(f_\theta(x))]\leq \alpha$ holds too. Since $t>0$, the restriction of the chance constraint can be alternatively expressed as $\mathbb{E}_\theta\left[[t+f_\theta(x)]_+\right]\leq \alpha t$. In fact, the requirement $t>0$ can be dropped because for all negative $t$, the constraint $\mathbb{E}_\theta[t+f_\theta(x)]_+\leq \alpha t$ becomes infeasible. And for $t=0$, the restricted constraint yields $\mathbb{E}_\theta\left[[f_\theta(x)]_+\right]\leq 0$ or equivalently $f_\theta(x)<0$ for all $\theta$, so that the original chance constraint holds trivially. Therefore, imposing the convex constraint $\mathbb{E}_\theta\left[[t+f_\theta(x)]_+\right]\leq \alpha t$ for some $t$ constitutes a restriction of the original chance constraint $\Pr[f_\theta(x)\geq 0]\leq \alpha$. 


By identifying $f_\theta(x)$ with $\underline{v}_n-v_n(\bq,\btheta)$ and introducing an auxiliary variable $\underline{t}_n$ per bus $n$, we can now restrict the voltage chance constraint in \eqref{eq:chancelow4} by imposing constraint
\[\mathbb{E}\left[[\underline{t}_n+\underline{v}_n-v_n(\bq,\btheta)]_+\right]\leq  \alpha\underline{t}_n,~\forall n\in\mathcal{N}.\]
The chance constraints of \eqref{eq:chancehigh4} on upper voltage limits can be treated similarly using variables $\overline{t}_n$'s. Collecting auxiliary variables $\{\underline{t}_n,\overline{t}_n\}_{n=1}^N$ in vectors $\underline{\bt}$ and $\overline{\bt}$ accordingly, we can now formulate the convex restriction of \eqref{eq:opf4} as
\begin{subequations}\label{eq:opf7}
\begin{align}
\min_{\bq\in\mcQ_{\btheta},\underline{\bt},\overline{\bt}}~&~\mathbb{E}[\ell(\bq,\btheta)]\\ 
\mathrm{s.to}~&~\mathbb{E}\left[[\underline{\bt}+ \underline{\bv}-\bv(\bq,\btheta)]_+- \alpha\underline{\bt}\right]\leq \bzero \label{eq:opf7ineq1}\\
~&~\mathbb{E}\left[[\overline{\bt}+ \bv(\bq,\btheta)-\overline{\bv}]_+ -\alpha\overline{\bt}\right]\leq \bzero.\label{eq:opf7ineq2}
\end{align}
\end{subequations}
For brevity, let the expressions inside the expectation operator of \eqref{eq:opf7ineq1}--\eqref{eq:opf7ineq2} be represented by $\underline{\bg}\left(\underline{\bt},\bv(\bq,\btheta)\right)$ and $\overline{\bg}\left(\overline{\bt},\bv(\bq,\btheta)\right)$. Similar to \eqref{eq:opf6}, problem \eqref{eq:opf7} can be tackled using stochastic primal-dual updates upon replacing the ensemble with sample averages over $K$ scenarios. Different from the averaged formulation however, the constraint functions $\underline{\bg}\left(\underline{\bt},\bv(\bq,\btheta)\right)$ and $\overline{\bg}\left(\overline{\bt},\bv(\bq,\btheta)\right)$ are non-differentiable with respect to $(\bv,\underline{\bt},\overline{\bt})$. We use their subgradients instead.
\begin{subequations}\label{eq:spdupdatchance}
\begin{align}
    \bw^{k+1}&:=\bw^{k}-\mu_w\left(\nabla_{\bw}\ell^k+ \left(\partial_{\bw}\underline{\bg}^k\right)^\top\underline{\blambda}^k+
    \left(\partial_{\bw}\overline{\bg}^k\right)^\top\overline{\blambda}^k\right)\label{eq:spdupdatechance:kappa}\\
    \underline{\bt}^{k+1}&:=\underline{\bt}^{k}-\mu_t\left(\partial_{\underline{\bt}}\underline{\bg}^k\right)^\top\underline{\blambda}^k\label{eq:spdupdatechance:xilow}\\
    \overline{\bt}^{k+1}&:=\overline{\bt}^{k}-\mu_t\left(\partial_{\overline{\bt}}\overline{\bg}^k\right)^\top\overline{\blambda}^k\label{eq:spdupdatechance:xihigh}\\
    \underline{\blambda}^{k+1}&:=\left[\underline{\blambda}^{k}+\mu_{{\lambda}}\underline{\bg}\left(\underline{\bt}^{k+1},\bv\left(\bq^{k+1},\btheta^{k}\right)\right)\right]_+\label{eq:spdupdatechance:dlow}\\
    \overline{\blambda}^{k+1}&:=\left[\overline{\blambda}^{k}+\mu_{{\lambda}}\overline{\bg}\left(\overline{\bt}^{k+1},\bv\left(\bq^{k+1},\btheta^{k}\right)\right)\right]_+.\label{eq:spdupdatechance:dup}
\end{align}
\end{subequations}

To compute the needed subgradients, recall that a subgradient of $f(x)=[x]_+$ can be found as
\[\partial f(x)= \left\{\begin{array}{ll}
0&,~x<0\\
1&,~x\geq 0
\end{array}\right.=\mathbbm{1}(x).\]
The subgradients involved in \eqref{eq:spdupdatechance:xilow}--\eqref{eq:spdupdatechance:xihigh} can be computed using the chain rule as
\begin{subequations}
\begin{align}
\partial_{\underline{\bt}}\underline{\bg}&=\diag\left(\mathbbm{1}\left(\underline{\bt}+\underline{\bv}-\bv\right)\right)-\alpha \bI_N\\
\partial_{\overline{\bt}}\overline{\bg}&=\diag\left(\mathbbm{1}\left(\overline{\bt}+\bv-\overline{\bv}\right)\right)-\alpha \bI_N
\end{align}
\end{subequations}
where the indicator functions here are applied entrywise and evaluate to vectors. In turn, the subgradients appearing in \eqref{eq:spdupdatechance:kappa} can be found as
\begin{subequations}
\begin{align}
  \partial_{\bw}\underline{\bg}&=\partial_{\bv}\underline{\bg}\cdot \nabla_{\bw}\bv\\ \partial_{\bw}\overline{\bg}&=\partial_{\bv}\overline{\bg}\cdot \nabla_{\bw}\bv.  
\end{align}
\end{subequations}

The Jacobian $\nabla_{\bw}\bv$ has already been computed in \eqref{eq:jacobvolt}, while the subgradients with respect to voltage magnitudes are
\begin{subequations}
\begin{align}
\partial_{\bv}\underline{\bg}&=-\diag\left(\mathbbm{1}\left(\underline{\bt}+\underline{\bv}-\bv\right)\right)\\
\partial_{\bv}\overline{\bg}&=\diag\left(\mathbbm{1}\left(\overline{\bt}+\bv-\overline{\bv}\right)\right).
\end{align}
\end{subequations}

\begin{remark}
During the training phase of the DNN, one may encounter ill-conditioned samples $\{\bphi^k,\btheta^k\}$ that, along with the DNN output $\bq^k$, cause the power flow solver to diverge. One can reduce the number of such ill-conditioned samples from the training data set by sampling close to the nominal benchmark values. Nonetheless, if such ill-conditioned samples do occur, a straightforward \emph{recourse} functionality could be adopted. First, the sample $\bphi^k$ is fed into the DNN to obtain $\bq^k$. Then, the PF solver is called for $\btheta^k$ and $\bq^k$. If the PF solver converges within a prescribed number of maximum iterations, we move on with the gradient calculations. Otherwise, we draw a new sample $\{\bphi^k,\btheta^k\}$ and repeat the process.
\end{remark}

\section{Gradient-Free Implementation}\label{sec:modelfree}
As discussed earlier, the primal-dual updates of Section \ref{sec:train} require computing the (sub)gradients of losses and voltages with respect to inverter reactive power injections. This section puts forth a gradient-free implementation of the DNN updates relying on a \emph{digital twin} of the feeder. This implementation does not require computing gradients, but only evaluating losses and voltages. The digital twin can be a power flow solver (GridLAB-D or OpenDSS), or a hardware emulator of the feeder. Once fed with all power injections (i.e., grid conditions $\btheta$ and inverter setpoints $\bq$), the digital twin returns the vector of nodal voltages $\bv$ and ohmic losses $\ell$.

Such gradient-free approach can be practically relevant under two settings. First, when one does not want to deal with gradients as calculating them is not straight-forward and can be prone to errors. Second, when the feeder model is complicated and computing gradients is quite complex. That could be the case in unbalanced multiphase grids; if grid devices (regulators, capacitors, transformer banks) are to be included; and/or when loads are represented by detailed ZIP or exponential models. In such cases, gradient calculations can be cumbersome. On other hand, reliable power flow modules do exist and can be used to evaluate the needed functions under either settings.

\begin{figure}[t]
    \centering
    \includegraphics[scale=0.37]{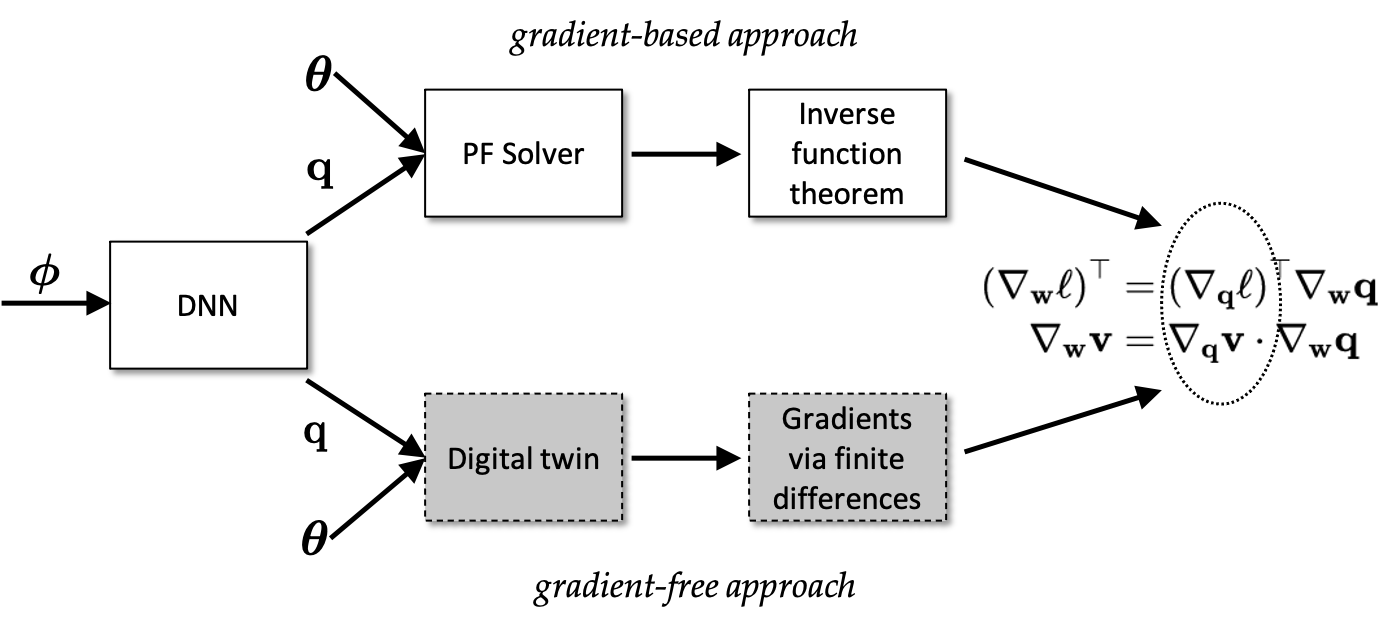}
    \caption{Steps for calculating $\nabla_{\bw}\ell$ and $\nabla_{\bw}\bv$ for the \emph{gradient-based} (top) and \emph{gradient-free} (bottom) approaches. The \emph{gradient-based} approach requires a solution to the power flow equations in conjunction with the inverse function theorem (Section \ref{sec:train}); the \emph{gradient-free} obtains the desired quantities by probing the digital twin and applying zero-order approximations.} 
    \label{fig:gf}
\end{figure}

To arrive at a gradient-free implementation, we cannot compute the partial derivatives appearing in the left-hand side of \eqref{eq:jacobvolt} and \eqref{eq:lossgrad}. Nonetheless, we can aim directly for the Jacobian matrix $\nabla_{\bq}\bv$ and the gradient vector $\nabla_{\bq}\ell$, and approximate them through finite differences. In detail, we resort to zeroth-order approximants (see~\cite{NesterovGF}) of the needed sensitivities by querying the digital twin twice to obtain two function evaluations as:
\begin{subequations}\label{eq:zerord}
\begin{align}
  \hat{\nabla}_{\bq}\ell&=\frac{\ell(\bq+\epsilon\cbq,\btheta)-\ell(\bq-\epsilon\cbq,\btheta)}{2\epsilon}\cbq^{\top}\label{eq:zerord1}\\
  \hat{\nabla}_{\bq} \bv&=\frac{\bv(\bq+\epsilon\cbq,\btheta)-\bv(\bq-\epsilon\cbq,,\btheta)}{2\epsilon}\cbq^{\top}\label{eq:zerord2}
\end{align}
\end{subequations}
where $\epsilon$ is the scale of perturbation, and $\cbq$ is a perturbation vector Gaussian distributed with zero-mean and standard deviation $\sigma_{\cbq}$. The quantities $\epsilon$ and $\sigma_{\cbq}$ are treated as hyper-parameters and are set during the training process.. The approximations in \eqref{eq:zerord} are carried out in three steps. First, the DNN is presented with the input $\bphi$ and its output $\bq$ is recorded. Second, the digital twin is presented with $(\bq,\btheta)$ and computes $\ell(\bq,\btheta)$ and $\bv(\bq,\btheta)$. Third, the digital twin is presented with $(\bq+\epsilon\cbq,\btheta)$ and computes $\bv(\bq+\epsilon\cbq,\btheta)$. Fig~\ref{fig:gf} compares the steps for obtaining the quantities $\nabla_{\bw}\ell$ and $\nabla_{\bw}\bv$ for the \emph{gradient-based} and \emph{gradient-free} approaches.

With the finite-difference approximants of \eqref{eq:zerord} in place, the \emph{gradient-free} primal-dual updates are straightforward to execute by approximating 
\begin{align*}
(\hat{\nabla}_{\bw}\ell)^\top&=(\hat{\nabla}_{\bq}\ell)^\top\nabla_{\bw}\bq\\
\hat{\nabla}_{\bw}{\bv}&=\hat{\nabla}_{\bq}{\bv}\cdot\nabla_{\bw}{\bq}
\end{align*}
where $\nabla_{\bw}{\bq}$ is again calculated gradient back-propagation. 
\section{Numerical Tests}\label{sec:tests}

\begin{figure}[t]
    \centering
    \includegraphics[scale=0.7]{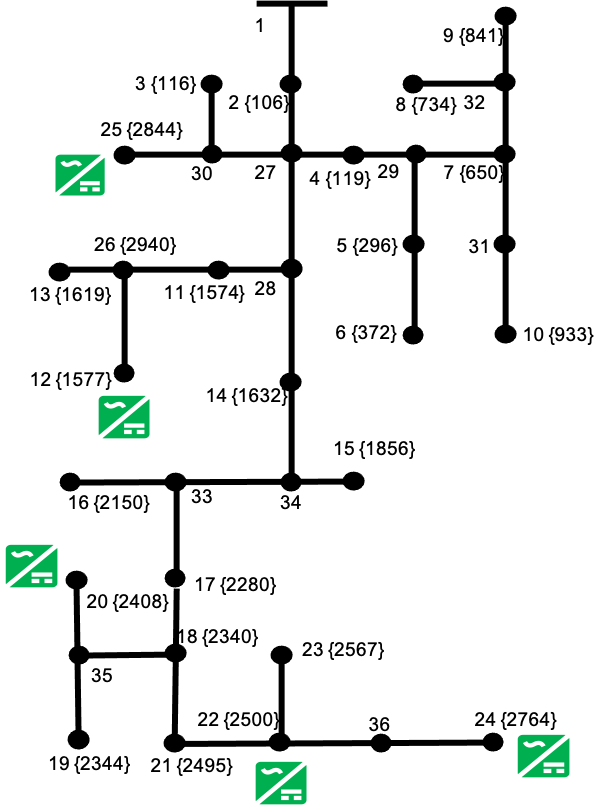}
    \caption{The IEEE 37-bus feeder used for the numerical tests. Node numbering follows the format \texttt{node number \{panel ID\}}. The inverters at nodes $\{12,20, 22, 24, 25\}$ provide reactive power control, whereas the rest operate at unit power factor.}
    \label{fig:ieee37}
\end{figure}

The performance of the proposed DNN-based control strategy was evaluated using a single-phase version of the IEEE 37-bus feeder. Real-world one-minute active load and solar generation data were extracted for April 2, 2011 from  the Smart* project~\cite{Smartmicrogrid},~\cite{Smartsolar}. For active loads, homes with IDs 20-369 were used. Averaged load demands were calculated by considering 10 homes at a time, and were serially allotted to buses 2-36 of Fig.~\ref{fig:ieee37}. The values of active loads were scaled so their maximum active load per node matched its benchmark value. Reactive loads were then added to each of these homes by sampling lagging power factors uniformly within $[0.9,1.0]$ and for each time interval. For solar generation, panel IDs were matched to the buses as shown in Fig.~\ref{fig:ieee37}. Solar generation values were scaled so the maximum generation per panel was $2$ times the benchmark value. Out of all the nodes with inverter-interfaced solar generation, those at nodes $\{12,20, 22, 24, 25\}$ were also providing reactive power support. The extracted data points were considered as available forecasts for $4$-hour control periods. Appropriate training and testing scenarios were created. Zero-mean white Gaussian noise was added to the $240$ one-minute data points from the forecast to create a total of $1200$ samples. The standard deviation of the Gaussian noise was set to $0.1$ times the mean load forecast. Out of the $1200$ samples, $960$ were used as the training set and the remaining $240$ formed the testing set. The training samples were additionally randomly shuffled to promote better generalization for the DNNs.

All tests were conducted on a $2.4$~GHz 8-Core Intel Core i9 processor laptop computer with $64$~GB RAM. Simulation scripts were written in Python and TensorFlow libraries to implement and train the DNNs. For the tests presented, four-layered fully connected DNNs were employed. The grid conditions vector $\btheta:=\left[\bp^g,\bp^c, \bq^c\right]^\top$ consisting of measurements from the $M\leq N $ buses equipped with smart meters were fed as inputs to the DNNs. Therefore, the input layers were chosen to have $3M$ neurons. The two subsequent hidden layers were fixed to having $3N$ and $2N$ neurons, respectively. Finally, the output layers had $5$ neurons corresponding to the $5$ inverters. All but the final layers of the DNNs employed the ReLU (rectified linear unit) activation with the final layers using a scaled \emph{tanh} activation to ensure the inverter limits $\bq^g\in\mathcal{Q}^t$. The weights for the DNN layers were initiated from a Gaussian distribution with zero mean and a unit standard deviation. The biases for the DNN layers, the dual variables, and the auxiliary variables were all initialized at zero. Additional modelling and training details are presented along with the discussions of the results as follows.

\subsection{Averaged Formulation}\label{subsec:averes}
For the average formulation, the DNN was fed with the complete vector of grid conditions $\btheta$ obtained from measurements collected at all buses. DNN weights were updated using {the DNN optimization algorithm Adam} with a learning rate of $0.001$. Dual variables were updated using SGD with a learning rate of $10$ that decayed with the square-root of the iteration index~\cite{LKMG17}. The model was then trained for $15$ epochs over the training scenarios.

\begin{figure}[t]
    \centering
    \includegraphics[scale=0.53]{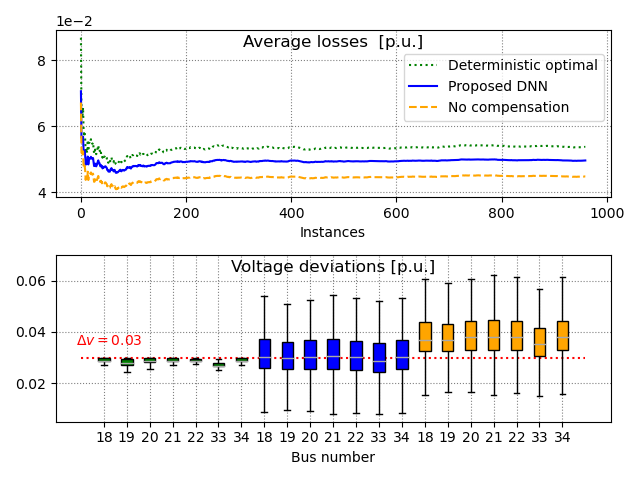}
    \caption{\emph{Top:} Time-averaged losses during the 12--4~pm training interval attained by the deterministic optimal control strategy of \eqref{eq:opf}; the proposed DNN-based inverter control; and no reactive power compensation by inverters. \emph{Bottom:} Box plots showing the first and third quantiles of the voltage deviations experienced across buses under the three control strategies. Due to high solar generation, the feeder experiences lower ohmic losses at the expense of severe over-voltages if there is no reactive power control by inverters. The deterministic optimal inverter control strategy regulates voltages by absorbing reactive power, which increases line currents and consequently losses. The proposed strategy achieves lower average losses over deterministic optimal inverter control as voltages are not constrained within $\pm3\%$ at all times.}
    \label{fig:avg_train}
\end{figure}

\begin{figure}[t]
    \centering
    \includegraphics[scale=0.53]{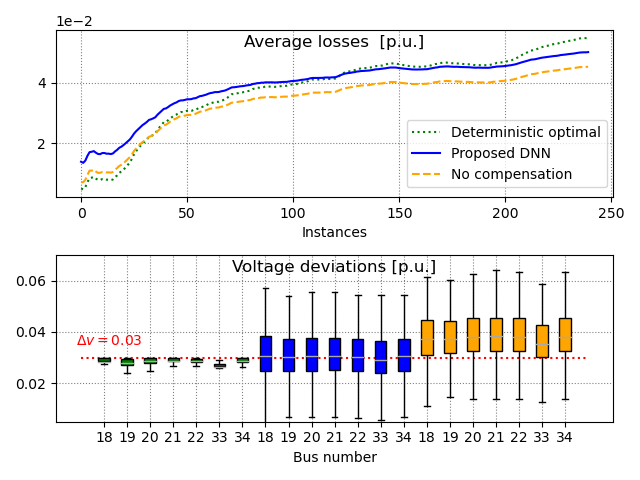}
   \caption{Results for averaged formulation over testing data during the interval 12--4~pm: Average losses under the deterministic optimal strategy, the proposed DNN-based approach; and no reactive power compensation are depicted on the top panel. Voltage deviations across buses under the three strategies are shown at the bottom panel.}
    \label{fig:avg_test}
\end{figure}

To demonstrate the efficacy of the proposed approach, the results are compared against a no-compensation scenario, i.e., the scenario where all inverters operate at unit power factor and provide no reactive power support. The DNN-based approach is also benchmarked against an deterministic optimal approach that solves the problem in \eqref{eq:opf} per minute. As discussed previously, such deterministic optimal approach might not be realistic to implement in real time due to the high computational burden. Fig.~\ref{fig:avg_train} compares the average losses and bus voltages under the three scenarios over the training set and during the high solar period of 12--4~pm. Without any reactive power compensation, buses $\{18,19,20,21,22,33,34\}$ experience over-voltages. The proposed DNN-based approach behaves as expected by lowering the average voltages at these buses down to the acceptable range. The deterministic optimal approach also achieves the same objective but by bringing all instantaneous voltages to the desired range whenever feasible. Note that both the DNN-based approach and the deterministic optimal approach incur higher losses when compared to the no compensation scenario. This is a result of increase in the magnitude of line currents on account of reactive power withdrawals. Since the deterministic optimal approach focuses on instantaneous voltage values rather than their averages, it incurs higher losses when compared to the DNN-based approach. The trained DNN was then evaluated over unseen scenarios of the testing set. As can be seen in Fig.~\ref{fig:avg_test}, the proposed approach performed remarkably well in maintaining voltages within limits and lowering average losses over the testing set.

To highlight the computational advantage of the DNN-based approach during operation, we compared it against the deterministic optimal strategy. The deterministic optimal strategy entails solving as many OPFs as the number of realization of $\btheta$'s, and was simulated by formulating a second-order conic program (SOCP) that models a convex relaxation of the original OPF~\cite{Low14}. The SOCP was solved for the $240$ samples of the testing dataset for $12-4$ pm in MATLAB using the SeDuMi solver. The simulation was found to take $171.24$~s to complete. On the other hand, the DNN-based approach, which requires only a forward pass through the DNN for each realization of $\bphi$, took only $0.85$~s to compute in Python. This indicates a significant computational speedup.

\begin{figure}[t]
    \centering
    \includegraphics[scale=0.37]{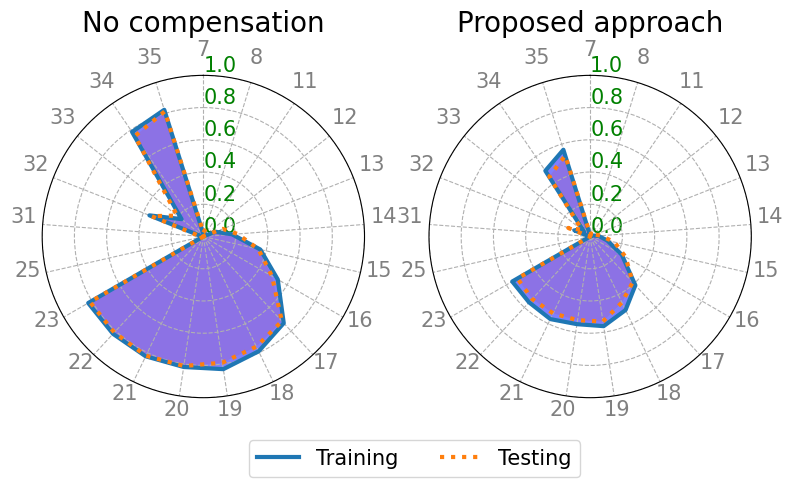}
    \caption{Voltage profiles during the 12--4~pm interval for averaged formulation. Results under no reactive power compensation and under the proposed DNN-based policies have been compared. Angular markings correspond to bus numbers, whereas radial markings are sampled probabilities of voltage violations.}
    \label{fig:polar_avg}
\end{figure}

Finally, the sample probabilities for voltage violations were calculated for the resulting voltages from the no compensation case and the proposed DNN-based approaches. For the buses with non-zero values for these probabilities, radar plots were drawn as shown in Fig.~\ref{fig:polar_avg}. The radial-axis represents the values for sample probabilities for voltage limit violations, whereas the angular markings correspond to the bus numbers. One can see that without any reactive power support, the grid faces a high probability of voltage violations over both training and testing. While the average formulation successfully regulates the average voltages, its effect on reducing the total number of occurrences of voltage violations is somewhat modest with buses $\{19,20,21,22,23\}$ violating the voltage limits for more than half of the scenarios.

\subsection{Probabilistic Formulation}
\begin{figure}[t]
    \centering
    \includegraphics[scale=0.27]{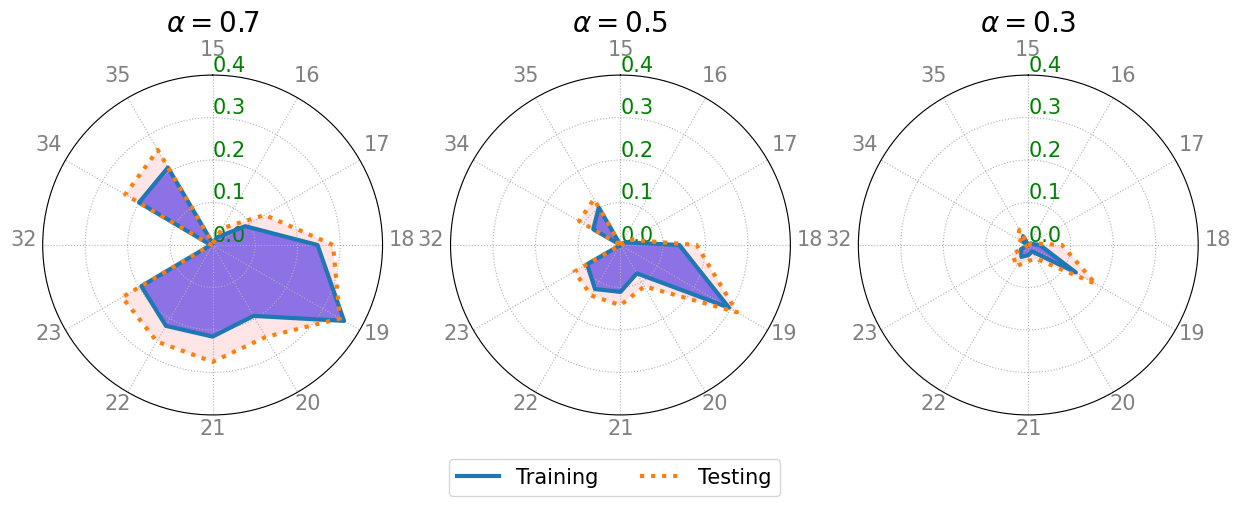}
    \caption{Results for probabilistic formulations for the 12--4~pm window. Voltage profiles for different values of $\alpha=\{0.7,0.5,0.3\}$ are depicted. Angular markings correspond to bus numbers whereas radial markings are sampled probabilities of voltage limits violation.}
    \label{fig:chance}
\end{figure}

Employing the DNN architecture from the previous subsection and the full input vector $\btheta$, updates in \eqref{eq:spdupdatchance} were applied to train the DNN for the probabilistic formulation. The DNN weights were updated using Adam with a learning rate of $0.001$. For the auxiliary variables $\{\underline{\bt},\overline{\bt}\}$, Adam optimizer with a learning rate of $0.001$ was deployed, and the dual variables were updated using SGD with a learning rate of $1$ that decayed with the square-root of the iteration index. The model required a higher number of $20$ epochs to converge during training because of the additional primal variables $\bt$. The experiments were conducted for the same time period of 12--4~pm. The experiments were repeated for three different values of $\alpha=\{0.7,0.5,0.3\}$ and the radar plots for the resulting sample probabilities of voltage violations are shown in Fig.~\ref{fig:chance}.

As desired, when compared to the averaged formulation, the occurrences of voltage violations under the probabilistic formulation were found to be drastically less for lower values of $\alpha$. Since the calculated sample probabilities came out to be less than the selected $\alpha$, the results in Fig.~\ref{fig:chance} confirm the conservative nature (being a convex restriction) of \eqref{eq:opf7}. 

\subsection{Gradient-Free Implementation}
\begin{figure}[t]
    \centering
    \includegraphics[scale=0.5]{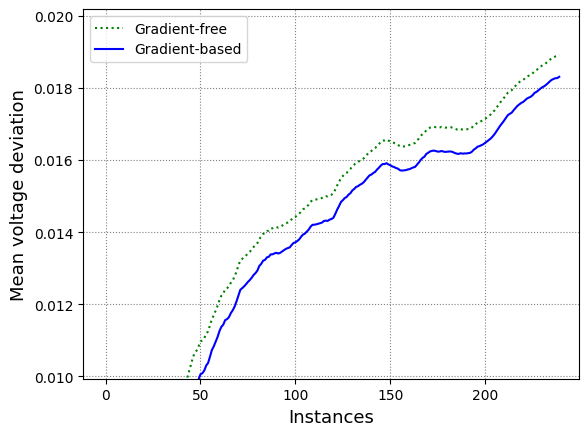}
    \caption{Mean voltage deviations across the buses under gradient-based and gradient-free approaches for the averaged formulation and during the testing 12--4~pm interval.}
    \label{fig:avg_mf}
\end{figure}

\begin{figure}[t]
    \centering
    \includegraphics[scale=0.4]{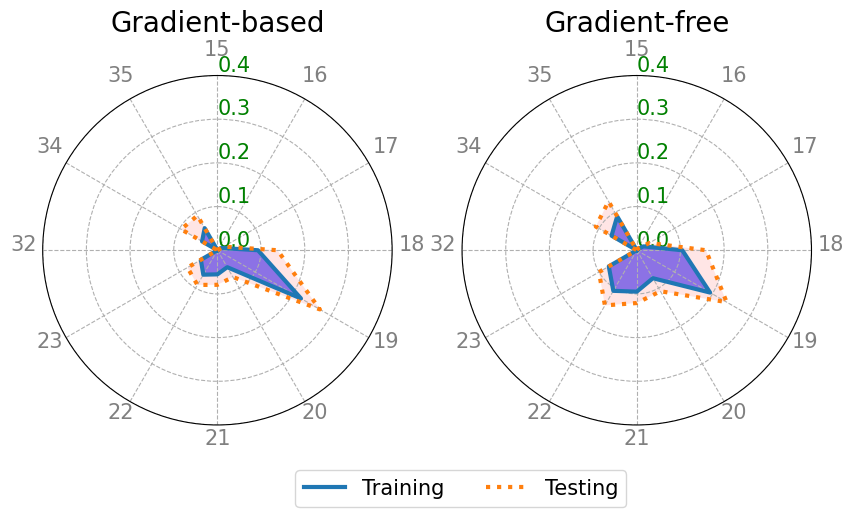}
    \caption{Comparison of gradient-based and gradient-free approaches for the probabilistic formulation with $\alpha=0.5$ over the testing 12--4~pm window. Voltage profiles for both approaches are depicted using polar plots.}
    \label{fig:chance_mf}
\end{figure}

To quantify the accuracy of the gradient approximations in \eqref{eq:zerord}, the DNNs for the averaged and probabilistic formulations were trained under the gradient-free fashion. The scale of perturbation $\epsilon$ was set to $0.1$, and the perturbations $\cbq$ were sampled from a zero-mean Gaussian distribution with a unit standard deviation. The gradient-free approaches were compared to their gradient-based counterparts over the same time periods and the results shown in Figs.~\ref{fig:avg_mf}-\ref{fig:chance_mf}. 

Fig.~\ref{fig:avg_mf} shows the mean voltage deviations across all the buses and time under the two approaches. The gradient-free approach incurs slightly higher voltage violations, but otherwise matches the performance of the gradient-based approach surprisingly well without any explicit knowledge of the underlying relationships. Similar results were confirmed for the probabilistic formulation in Fig.~\ref{fig:chance_mf}, where both approaches yield similar sample probabilities for voltage violations over training and testing, with $\alpha=0.5$.

\subsection{Partial Inputs}
\begin{figure}[t]
    \centering
    \includegraphics[scale=0.55]{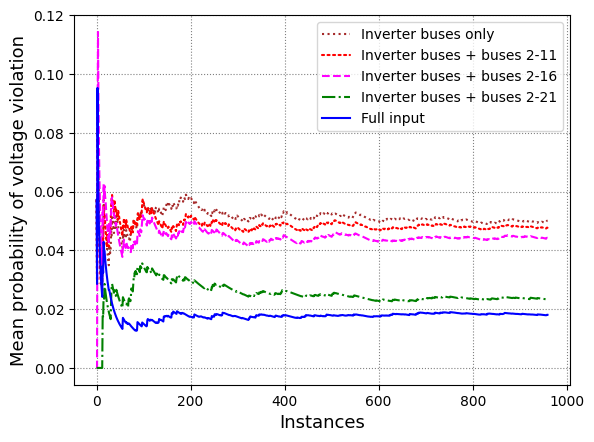}
    \caption{Mean sampled probabilities of voltage violations for the probabilistic formulation with $\alpha=0.5$ during training over 12--4~pm. The dimension of the DNN input increases as the number of nodes monitored in real time increases.}
    \label{fig:partial_tr}
\end{figure}

\begin{figure}[t]
    \centering
    \includegraphics[scale=0.55]{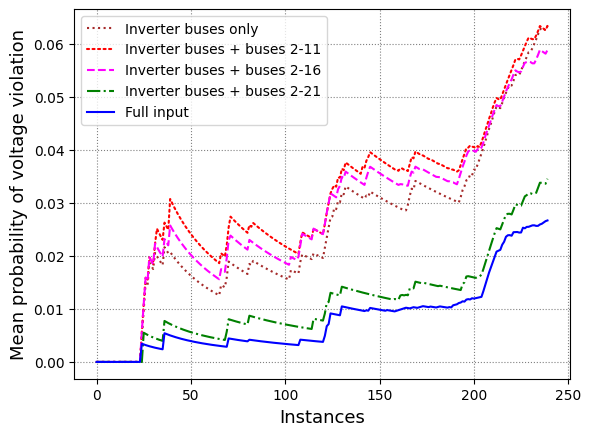}
    \caption{Mean sampled probabilities of voltage violations for the probabilistic formulation with $\alpha=0.5$ during testing over 12--4~pm. The dimension of the DNN input increases as the number of nodes monitored in real time increases.}
    \label{fig:partial_test}
\end{figure}

To study the effect of partial DNN inputs on the control performance of the learned DNN policy, we varied the number of nodes whose data are telemetered in real time for the probabilistic formulation. First, real-time meters were assigned to all nodes with solar generation. Then, different input scenarios were simulated by expanding the subset of inverter-equipped nodes with real-time metering moving from nodes 2--11; nodes 2--16; nodes 2--21; and a full input vector consisting of all nodes. The mean-sampled probabilities of voltage violations across time and buses were recorded. The value of $\alpha$ was set to $0.5$ for all scenarios. Figures~\ref{fig:partial_tr} and \ref{fig:partial_test} show the results recorded respectively during training and testing. The results confirm the intuition that the performance of the control policies improves as more real-time information is provided to the DNN. As the nodes with real-time monitoring increase, the sampled probabilities of voltage violations decrease.

\section{Conclusions}\label{sec:conclusions}
This work has presented a DNN-based approach for stochastic optimal inverter control. To capture uncertainty, the grid conditions have been modeled as random variables and the associated inverter setpoints by a stochastic policy learned through a DNN. The DNN is periodically trained offline to minimize the average ohmic losses and maintain either the average voltages within limits or the probability of voltage deviation occurrences low. Training is accomplished by adopting existing stochastic primal-dual updates and their gradient-free counterparts to the AC OPF setup. The proposed scheme not only expedites the computation of near-optimal inverter setpoints, but also resolves two practical difficulties: \emph{d1)} How to solve an OPF using only a power flow solver or a digital twin of the feeder? and \emph{d2)} How to deal with an OPF if grid conditions are only partially known? 

Numerical tests on the benchmark IEEE 37-bus feeder showcase the salient features of the novel methodology. The adopted stochastic primal-dual updates train a DNN-based policy using a modest number of training samples. The DNN is numerically shown to produce an inverter dispatch that minimizes the OPF cost while satisfying the operating constraints. Although the averaged formulation succeeds in maintaining the voltage at each node within limits on the average, violations do occur frequently. To that end, the chance-constrained formulation may be more relevant. The latter comes at minimal extra complexity over the averaged formulation. Being a convex restriction the numerical tests corroborate that it is a conservative yet safe scheme. Our experiments have also shown how the DNN-based policy can be driven with incomplete grid conditions, and demonstrated the improvement in feasibility when more information is presented to the DNN. A final interesting outcome is that when the DNNs are trained using the gradient-free updates, the degradation in performance is minimal although the learner has less information about the feeder at its disposal. 

Some open questions are to: \emph{i)} Extend the implementation to a multiphase grid model with detailed ZIP loads, regulators, and capacitor banks; \emph{ii)} Consider additional network constraints and/or alternative objectives; \emph{iii)} Experiment with the frequency of re-training the DNN, the size of the training dataset, and the way it has been generated; \emph{iv)} Utilize the DNN output only to warm-start an actual OPF solver; and \emph{v)} Optimally select the grid proxies to improve inverter control performance under a communication budget.


\bibliographystyle{IEEEtran}
\bibliography{myabrv,inverters}
\end{document}